\documentclass[12pt]{amsart}

\usepackage{amsmath,amssymb}
\usepackage{amsfonts}
\usepackage{amsthm}
\usepackage{latexsym}
\usepackage{graphicx}
\usepackage{epsfig}

\makeatletter \@addtoreset{equation}{section} \makeatother

\def\ddt{\frac{d}{dt}}
\def\ppt{\frac{\partial}{\partial t}}
\numberwithin{equation}{section}

\def\bg{{\mathbf g}}

\newtheorem{prop}{Proposition}[section]
\newtheorem{theo}[prop]{Theorem}

\newtheorem{coro}[prop]{Corollary}
\newtheorem{rema}[prop]{Remark}
\def\Dint{\displaystyle\int}

\begin{document}

\title{First Eigenvalues of Geometric Operators under the Ricci Flow}

\author{Xiaodong Cao}
\thanks{Research partially supported by an MSRI postdoctoral fellowship }

\address{Department of Mathematics,
  Cornell University, Ithaca, NY 14853-4201}
\email{cao@math.cornell.edu}

\renewcommand{\subjclassname}{%
  \textup{2000} Mathematics Subject Classification}

\date{Oct. 5th,  2007}

\maketitle

\markboth{Xiaodong Cao} {First Eigenvalues of Geometric Operators
under the Ricci Flow }

\begin{abstract}  In this paper, we prove that the first
eigenvalues of  $-\Delta + cR$ ($c\geq \frac14$) is nondecreasing
under the Ricci flow. We also prove the monotonicity under the
normalized Ricci flow for the case $c= 1/4$ and $r\le 0$.
\end{abstract}

\section{\textnormal{\bf First Eigenvalue of $-\Delta + cR$ }}

Let $M$ be a closed Riemannian manifold, and $(M,\bg(t))$ be a
smooth
 solution to the Ricci flow equation $$\ppt \bg_{ij}=-2R_{ij}$$ on $0\le t<T$.
 In \cite{C},
 we prove that all eigenvalues $\lambda (t)$  of the operator
 $-\Delta + \frac{R}{2}$ are nondecreasing under the Ricci flow on manifolds with
 nonnegative curvature operator. Assume $f=f(x,t)$ is the
 corresponding eigenfunction of
 $\lambda (t)$, that is
 $$(-\Delta + \frac{R}{2})f(x,t)=\lambda (t) f(x,t)$$ and
  $\Dint_M f^2d\mu = 1$.  More generally,  we define
\begin{equation}\label{defn1}
\lambda(f,t)=\Dint_M (-\Delta f + \frac{R}{2}f)f d\mu,
\end{equation}
 where $f$ is
a smooth function satisfying $$\frac{d}{dt}(\Dint_M f^2 d\mu)=0,~
\Dint_M f^2 d\mu=1.$$  We can then derive the monotonicity formula
under the Ricci flow.

\begin{theo}\label{thm1}\cite{C}
On a closed Riemannian manifold with nonnegative curvature
operator, the eigenvalues of the operator $-\Delta + \frac{R}{2}$
are nondecreasing under the Ricci flow. In particular,
\begin{equation} \label{equ1}
\textstyle{\frac{d}{dt}\lambda(f,t) = 2\Dint_M R_{ij}f_if_j d\mu +
\Dint_M |Rc|^2f^2d\mu\ge 0.}
\end{equation}
\end{theo}

In  (\ref{equ1}), when $\frac{d}{dt}\lambda(f,t)$ is evaluated at
time $t$, $f$ is the corresponding eigenfunction of $\lambda(t)$.

\begin{rema}\label{rema1}
Clearly, at time $t$, if $f$ is the eigenfunction of the
eigenvalue $\lambda(t)$, then $\lambda(f,t) = \lambda(t)$. By the
eigenvalue perturbation theory, we may assume that there is a
$C^1$-family of smooth eigenvalues and eigenfunctions (for
example, see \cite{kleinerlott}, \cite{rs} and \cite{chowetc1}).
 When $\lambda$ is the lowest eigenvalue, we can further
assume that the corresponding eigenfunction $f$ be positive. Since
the above formula does not depend on the particular evolution of
$f$, so $\ddt \lambda(t)=\ddt \lambda(f,t)$.
\end{rema}

\begin{rema}\label{rema2}
In \cite{ljf1}, J. Li used the same technique to prove that the
monotonicity of the first eigenvalue of $-\Delta + \frac12 R$
under the Ricci flow without assuming nonnegative curvature
operator. A similar result appeared in the physics literature
\cite{osw05}.
\end{rema}

\begin{rema}\label{rema3}
When $c=\frac14$, the monotonicity of first eigenvalue has been
established by G. Perelman in \cite{perelman1}. The evolution of
first eigenvalue of Laplace operator under the Ricci flow has been
studied by L. Ma in \cite{ma1}. The evolution of Yamabe constant
under the Ricci flow has been studied by S.C. Chang and P. Lu in
\cite{lu1}.
\end{rema}

In this paper, we shall study the first eigenvalues of operators
$-\Delta + cR$ ($c\geq \frac14$) without {\em curvature
assumption} on the manifold. Our first result is the following
theorem.

\begin{theo}\label{thm2}
Let ($M\sp n,\bg(t)$), $t\in[0,T)$, be a solution of the Ricci
flow on a closed Riemannian manifold $M\sp n$. Assume that
$\lambda (t)$ is the lowest eigenvalue of  $-\Delta + c R$ ($c\geq
\frac14$),
 $f=f(x,t)>0$
satisfies
\begin{equation}\label{equ2}
-\Delta f(x,t) + c R f(x,t)=\lambda(t)f(x,t)
\end{equation}  and
$\textstyle{\Dint_M f^2(x,t) d\mu = 1}$. Then under the Ricci
flow, we have
\begin{equation} \nonumber \label{equ11}
\begin{array}{rll}
\textstyle{\frac{d}{dt}\lambda(t) = \frac12 \Dint_M
|R_{ij}+\nabla_i\nabla_j \varphi|^2\ e^{-\varphi}d\mu +
\frac{4c-1}{2} \Dint_M |Rc|^2\ e^{-\varphi}d\mu \ge 0,}
\end{array}
\end{equation}
where $e^{-\varphi} = f^2$.
\end{theo}

\begin{proof}({\bf Theorem \ref{thm2}})
Let $\varphi$ be a function satisfying $ e^{-\varphi(x)}=f^2(x)$.
We proceed as in \cite{C}, we have
\begin{equation}\label{equ12}
\begin{array}{rll}
\frac{d}{dt}\lambda(t)=&(2c-\frac12) \Dint R_{ij}\nabla_i
\varphi\nabla_j\varphi e^{-\varphi}d\mu\\ & -(2c-1)\Dint R_{ij}
\nabla_i \nabla_j \varphi e^{-\varphi} d\mu + 2c\Dint
|Rc|^2e^{-\varphi}d\mu.
\end{array}
\end{equation}
Integrating by parts and applying the Ricci formula, it follows
that
\begin{equation}\label{equ14}
\Dint R_{ij}\nabla_i \nabla_j\varphi\ e^{-\varphi}d\mu = \Dint
R_{ij}\nabla_i \varphi\nabla_j\varphi\ e^{-\varphi}d\mu
-\frac{1}{2}\Dint R\Delta e^{-\varphi}d\mu\\
\end{equation}
and
\begin{equation}\label{equ8}
\begin{array}{rll}
\Dint R_{ij}\nabla_i\nabla_j \varphi e^{-\varphi}d\mu &+ \Dint
|\nabla\nabla \varphi|^2 e^{-\varphi}d\mu \\
= &-\Dint\Delta e^{-\varphi}\big( \Delta\varphi+\frac{1}{2}R
-\frac{1}{2}|\nabla\varphi|^2 \big)d\mu\\
= &(2c-\frac12)\Dint R \Delta e^{-\varphi}d\mu.\\
\end{array}
\end{equation}
In the last step, we use
\begin{equation}\nonumber
2\lambda(t) = \Delta\varphi+2c R -\frac{1}{2}|\nabla\varphi|^2.
\end{equation}
Combining (\ref{equ14}) and (\ref{equ8}), we arrive at
\begin{equation}\label{eqn10}
\Dint |\nabla\nabla \varphi|^2 e^{-\varphi}d\mu= 2c \Dint R \Delta
e^{-\varphi}d\mu - \Dint R_{ij}\nabla_i \varphi\nabla_j\varphi\
e^{-\varphi}d\mu.
\end{equation}
Plugging (\ref{eqn10}) into (\ref{equ12}), we have
\begin{equation}\nonumber
\begin{array}{rll}
\frac{d}{dt}\lambda(t)=& \Dint R_{ij} \nabla_i \nabla_j \varphi
e^{-\varphi} d\mu + 2c\Dint |Rc|^2e^{-\varphi} d\mu \\
&+c\Dint R\triangle (e^{-\varphi}) d\mu -\frac12 \Dint
R_{ij}\nabla_i \varphi\nabla_j\varphi\
e^{-\varphi}d\mu\\
=& \Dint R_{ij}\nabla_i \nabla_j\varphi e^{-\varphi}d\mu + 2c
\Dint |Rc|^2e^{-\varphi}d\mu + \frac12\Dint
|\nabla \nabla  \varphi|^2 e^{-\varphi}d\mu\\
=& \frac12 \Dint |R_{ij}+\nabla_i\nabla_j \varphi|^2\
e^{-\varphi}d\mu + (2c-\frac12) \Dint |Rc|^2\ e^{-\varphi}d\mu \ge
0.
\end{array}
\end{equation}
This proves the theorem as desired.
\end{proof}

\section{\textnormal{\bf  First Eigenvalue under the Normalized Ricci Flow}}

In this section, we derive the evolution of $\lambda (t)$ under
the normalized Ricci flow equation $$\ppt
\bg_{ij}=-2R_{ij}+\frac{2}{n} r\bg_{ij}.$$ Here
$\textstyle{r=\frac{\int_M R d\mu}{\int_M d\mu}}$ is the average
scalar curvature. It follows from Eq. (\ref{equ2}) that $\lambda
\le c r$. We now compute the derivative of the lowest eigenvalue
of $-\Delta + cR$.
\begin{theo}\label{thm3}
Let ($M\sp n,\bg(t)$), $t\in[0,T)$, be a solution of the
normalized Ricci flow on a closed Riemannian manifold $M\sp n$.
Assume that $\lambda (t)$ is the lowest eigenvalue of  $-\Delta +
c R$ ($c\geq \frac14$), $f>0$ is the corresponding eigenfunction.
Then under the normalized Ricci flow, we have
\begin{equation} \label{equ21}
\begin{array}{rll}
\frac{d}{dt}\lambda(t)=-\frac{2r\lambda}{n}+\frac12 \Dint_M
|R_{ij}+\nabla_i\nabla_j \varphi|^2\ e^{-\varphi}d\mu +
\textstyle{\frac{4c-1}{2}} \Dint_M |Rc|^2\ e^{-\varphi}d\mu ,
\end{array}
\end{equation}
where $e^{-\varphi} = f^2$. Furthermore, if
$\textstyle{c=\frac14}$ and $r\le 0$, then
\begin{equation} \label{equ22}
\begin{array}{rll}
\frac{d}{dt}\lambda(t) =\frac{2}{n}r(\textstyle{\lambda
-\frac{r}{4}})+\frac12 \Dint_M |R_{ij}+\nabla_i\nabla_j
\varphi-\frac{r}{n}g_{ij}|^2\ e^{-\varphi}d\mu \ge 0.
\end{array}
\end{equation}
\end{theo}

\begin{rema}After we submitted our paper,
J. Li suggested to us that (\ref{equ22}) is true for all $c\geq
1/4$, with an additional nonnegative term $$
\textstyle{\frac{4c-1}{2} \Dint_M |Rc-\frac{r}{n}g_{ij}|^2\
e^{-\varphi}d\mu },$$ see \cite{ljf2} for a similar result.
\end{rema}

\begin{rema}\label{rema2.2}
As a consequence of the above monotonicity formula of
$\lambda(t)$, we can prove that both compact steady and expanding
Ricci breathers (cf. \cite{Isoliton}, \cite{perelman1}) must be
trivial, such results have been discussed by many authors (for
example, see \cite{Isoliton}, \cite{Hsurvey}, \cite{Hsurface},
\cite{perelman1}, \cite{C} and \cite{ljf1}, etc.).
\end{rema}

When $M$ is a two-dimensional surface, $r$ is a constant. We have
the following corollary.

\begin{coro}
Let ($M\sp 2,\bg(t)$), $t\in[0,T)$, be a solution of the
normalized Ricci flow on a closed Riemannian surface $M^2$. Assume
that $\lambda (t)$ is the lowest eigenvalue of $-\Delta + c R$
($c\geq \frac14$), we have $e^{rt}\lambda$ is nondecreasing under
the normalized Ricci flow. Moreover, if $r\le 0$, then $\lambda$
is nondecreasing.
\end{coro}

{\bf Acknowledgement:} The author would like to thank Dr. Junfang
Li for bringing \cite{ljf1} and \cite{ljf2} to his attention and
for helpful discussion on this subject. He would like to thank
Professor Eric Woolgar for his interest and pointing out the
reference \cite{osw05}. He would like to thank Professor Duong H.
Phong for his interest and encouragement. He also wants to thank
MSRI for the generous support.

\bibliographystyle{halpha}
\bibliography{bio}
\end{document}